\documentclass[]{statsoc}
\usepackage[latin1]{inputenc}
\usepackage{amsmath}
\usepackage{amsfonts}
\usepackage{amssymb}
\usepackage{mathrsfs}
\usepackage{enumerate}

\def\be{\begin{equation}}
\def\ee{\end{equation}}
\def\bea{\begin{eqnarray}}
\def\eea{\end{eqnarray}}

\def\Prob{{\rm Prob}}
\def\tF{{\tau_{\scriptscriptstyle F}}}
\def\tG{{\tau_{\scriptscriptstyle G}}}

\title[Combining independent $P$-values]{Combining independent, arbitrarily weighted $P$-values: 
a new solution to an old problem using a novel expansion with controllable accuracy}
\author{Gelio Alves}
\address{National Institutes of Health, Bethesda, USA\\ [3pt]\hspace*{1in}and}

\author[Alves and Yu]{Yi-Kuo Yu$^\dag$} 
 
\address{National Institutes of Health, Bethesda, USA}  

\begin{document}
%\date{August 24th, 2010}
\begin{abstract}
Good's formula and Fisher's method are frequently used for combining independent $P$-values. 
  Interestingly, the equivalent of Good's formula already emerged in 1910
 and mathematical expressions relevant to even more general situations have been 
 repeatedly derived, albeit in different context. 
  We provide here a novel derivation and show how the analytic formula obtained
  reduces to the two aforementioned ones as special cases. 
  The main novelty of this paper, however, is the explicit treatment of 
  nearly degenerate weights, which are known to cause numerical instabilities. 
  We derive a controlled expansion, in powers of differences in inverse weights,     
  that provides both accurate statistics and stable numerics.
\end{abstract}

\maketitle

\footnotetext{\\ {\it $\dag$Address for correspondence}: National Center for
 Biotechnology Information, National Library of Medicine, National Institutes of Health, 
 8600 Rockville Pike, Bethesda, MD 20894, USA. \\
 E-mail: yyu@ncbi.nlm.nih.gov}

\section{Introduction}
The question of how to obtain an overall significance level for the results 
of {\it independent} runs of studies has been investigated since the 1930s~\citep{Tippett_1931,Fisher_1932,Pearson_1933,Pearson_1938}. 
In fact, forming a single statistical significance out of multiple independent tests
 has been an important subject of study in numerous area of scientific disciplines, including 
  social psychology~\citep{Stouffer_1949,MB_1954}, medical research~\citep{Olkin_1995}, genetics~\citep{LDGB_2001}, proteomics~\citep{AWWSetal_2008}, genomics~\citep{HI_2007}, bioinformatics~\citep{BG_1998,YGASetall_2006} and others.

Frequently used methods for combining independent $P$-values fall into
numeric and analytic categories. This classification is not totally precise
 since method such as Fisher's started out with the necessity of inverting the $\chi^2$ 
  cumulative distribution function and thus seemed like a numerical approach~\citep{Pearson_1938}. 
The method mentioned in~\citep{Bharucha-Reid_1960}, although not in the context of 
 combining $P$-values, brought out an analytical expression for combined $P$-value using Fisher's method, 
   thus effectively brought Fisher's method into analytic category.  
   In the context of combining $P$-values, by mapping to a known result 
   by Feller~\citep{Feller_1966},~\citet{BG_1998} also 
  provided an analytic formula for Fisher's combined $P$-value. 
 Numerical approaches typically involve inverting cumulative distribution functions. 
 For example, Stouffer's z-methods~\citep{Stouffer_1949}, whether 
 unweighted~\citep{Whitlock_2005} or weighted~\citep{Liptak_1958,KT_1994}, require inverting the error function. 
  Lancaster's generalization~\citep{Lancaster_1961,Koziol_1996} of Fisher's formalism also requires inverting 
   gamma distribution function to incorporate unequal weighting for $P$-values combined.   

 In this paper, we focus on analytic methods only. Two existing analytic approaches, 
 Fisher's~\citep{Fisher_1932,Bharucha-Reid_1960,BG_1998,AWWSetal_2008} and Good's~\citep{Good_1955,Likes_1967} , are frequently employed. 
  Fisher's method combines $L$ independent tail-area probabilities
  democratically to form a single significance assignment while Good's formula weights 
  each tail-area probability differently to form a different single significance assignment.
 Since Good's expression involves, in the denominator, pairwise differences between weights, he
  cautiously remarked that the expression may become ill-conditioned and thus the calculations
   should be done by holding more decimal places when weights of similar magnitudes exist. 
  This statement has been paraphrased by many authors~\citep{Bhoj_1992,OS_2001,Hou_2005}. 
    
In addition to the cases considered by Fisher and Good, it is foreseeable that one may wish to
 categorize obtained independent $P$-values into different groups so that one would like
  to weight $P$-values within the same group democratically and weight different group differently.   
 We will call this scenario the {\it general case} (GC). The GC naturally occurs since one may wish      
  to categorize data obtained from the same type of experimental instruments into the same group,
   and data collected from different instrument types may justify the use of different weights. 
 When there is only one instrument type, the GC reduces to Fisher's consideration.
  When there exist no replicates within each instrument types, the GC coincides 
   with the consideration of Good. In principle, the weighted version of the Stouffer's method can also
    be used for this purpose. Since the main scope here is to pursue analytic results, we
  won't delve into methods in numerical category.    
            
It is important, however, to point out that the mathematical problem of combining $P$-values is also
 related to other areas of research. For example, the equivalent of Good's formula 
 had already emerged in 1910 in the context of sequential radioactive decay~\citep{Bateman_1910},
  while analytical expression for Fisher's combined $P$-value emerged in 1960 as a special case of the
   former when all the decay constants are identical~\citep{Bharucha-Reid_1960}. 
 After Good's formula~\citep{Good_1955}, the same distribution function was rederived by 
 ~\citet{MG_1965} and later on by~\citet{Likes_1967}.     
 As for the GC, Fisher's method included, the mathematical equivalents appear in different areas of 
 studies mainly under the consideration of sum of exponential/gamma variables. 
 The distribution function of linear combination of exponential/gamma variables 
 are useful in various fields. When limited to 
  exponential variables, it results in the Erlang distribution that is often encountered
   in queuing theory~\citep{Morse_1958}. It is also connected to the renewal theory~\citep{Cox_1962}, 
   time series problem~\citep{MacNeill_1974}, and can be applied to model reliability~\citep{JK_2003}. 
   The intimate connections between these seemingly different problems are not obvious at first glance. 
   Consequently, it is not surprising that
  the distribution function has been rediscovered/rederived many times and that some information about it
   has not been widely circulated. Our literature searches show that the first explicit 
   result (without further derivatives involved) for the distribution function was obtained by \citet{Mathai_1983}. 
  Subsequently,  motivated by different contexts,~\citet{Harrison_1990},~\citet{AM_1997}, and ~\citet{JK_2003} 
  all rederived the same distribution function. 
 
Employing a complex variable integral formulation, we are able to provide a different derivation 
 of the distribution function and become the first to make connection to the GC of combining $P$-values.  
%  an exact significance formula for the general case described.
 Since both Fisher's and Good's considerations arise as special limiting cases of the GC, 
 we also illustrate that our cumulative probability distribution for GC indeed reduces to the appropriate limiting 
   formulas upon taking appropriate parameters. 
 The main novelty of this paper, however, is the explicit treatment of cases
   where nearly degenerate weights exist. These cases are known for numerical instabilities, which 
    motivated many authors to pursue uncontrolled approximations~\citep{HS_1977,Gabler_1987,Bhoj_1992,Hou_2005}.     
   We have derived a controlled expansion, in power of differences in inverse weights,     
  that provides both accurate statistics and stable numerics.

In the following section, we will first summarize Fisher's and Good's methods for combining $P$-values, followed by
 the mathematical definition of the GC. A section devoted to derivation of the probability distribution function
  and cumulative probability for the GC then follows. We then delve into the case of nearly degenerate weights 
  and provides a formula with controllable accuracy for combining $P$-values. 
 A few examples of using the main results are then provided. 
This paper concludes with future directions.  
  A C++ program, which computes the combined $P$-values with equal numerical stability 
  regardless of whether weights are (nearly) degenerate or not, 
 is available upon request from the authors.

\section{Summary of Fisher's and Good's methods for combining $P$-values}
Imagine that one wishes to combine $L$ independent $P$-values $p_1, p_2, \ldots, p_L$, 
 each of which is drawn from an uniform, independent distribution over $(0,1]$. 
 For later convenience, let us define 
\bea  
 \tF &\equiv&  p_1\cdot p_2 \cdots p_L  \; ,\label{Fisher_tau} \\
 \tG &\equiv& p_1^{w_1} \cdot p_2^{w_2} \cdots p_L^{w_L} \label{Good_tau} \;.  
\eea 
 To form  a unified significance, Fisher and Good considered respectively 
  the stochastic quantities $Q_F$ and $Q_G$, defined by
\bea 
 Q_F &\equiv&  x_1 \cdot x_2 \cdots x_L \label{Fisher_Q} \; , \\
 Q_G &\equiv& x_1^{w_1} \cdot x_2^{w_2} \cdots x_L^{w_L} \label{Good_Q} \; , 
\eea 
where each $x_i$ represents a random variable drawn from an uniform, independent distribution over $(0,1]$.
 The following probabilities 
\bea
\Prob (Q_F \le \tF) &=& \tF \sum_{l=0}^{L-1} {[\; \ln (1/\tF)\,]^{\,l} \over l! } \label{Comb.Fisher}\\
\Prob (Q_G \le \tG) &=&  \sum_{l=1}^L\,  \Lambda_l \; \tG^{1/w_l} \label{Comb.Good}
\eea
provide the unified statistical significances, corresponding respectively to Fisher's and Good's considerations,
 from combining $L$ independent $P$-values. 
 In eq.~(\ref{Comb.Good}), the prefactor $\Lambda_l$ is given by
\be  
\Lambda_l = \frac{w_l^{L-1}}{\prod_{k\ne l}(w_l - w_k)}.
\ee
Apparently, $\Lambda_l$ is ill-defined when the weight $w_l$ coincides with or is 
numerically close to any other weights $w_k$. 
 Although Fisher did not derive (\ref{Comb.Fisher}), from this point on, we shall refer to (\ref{Comb.Fisher})
  as Fisher's formula and (\ref{Comb.Good}) as the Good's formula.

\section{General case including Fisher's and Good's formulas}  
 Let us divide the $L$ $P$-values into $m$ groups with $1\le m \le L$.
 Within each group $k$, we weight the $n_k$ $P$-values within equally;
  while $P$-values in different groups are weighted differently. Therefore, when 
  $m=L$ and $n_k=1$ $\forall~k$, we have the Good's case and when $m=1$ and $n_1 = L$, we 
   reach Fisher case. We will therefore define the following quantities of interest
 \bea 
\tau &\equiv & \prod_{k=1}^m \left[ \prod_{j=1}^{n_k} p_{k;j} \right]^{w_k} \label{Unify_tau} \\
Q &\equiv &  \prod_{k=1}^m \left[ \prod_{j=1}^{n_k} x_{k;j} \right]^{w_k} \label{Unify_tau}
 \eea   
where each $x_{k;j}$ represents again a random variable drawn from an uniform, independent distribution over $(0,1]$.   
 The quantity of interest $\Prob(Q \le \tau)$, if obtained, should cover both results of 
  Fisher's and Good's as the limiting cases.  In the next section, we will derive
 an exact expression for $\Prob (Q \le \tau)$ and describe how to recover the results of 
 Fisher's and Good's. 
   
\section{Derivation of $\Prob (Q \le \tau)$}  
\label{sec:derivation}  
Let $F(\tau) \equiv \Prob (Q \le \tau)$, we then may write
\be
F(\tau) = \int_{0}^{1} \cdots  \int_{0}^{1}   \theta \left( \tau - \prod_{k=1}^m \left[ \prod_{i=1}^{n_k} x_{k;i} \right]^{w_k}\right) 
 \prod_{k=1}^m \prod_{j=1}^{n_k} dx_{k;j} \; , \label{F.def}
\ee
where the heaviside step function $\theta (x)$ takes value $1$ when $ x > 0$ and value $0$ when $x < 0$.    
Upon taking a derivative with respect to $\tau$, we obtain
\be
f(\tau) \equiv \frac{dF(\tau)}{d\tau} = 
 \int_{0}^{1} \cdots \int_{0}^{1}  \delta \left( \tau - \prod_{k=1}^m \left[ \prod_{j=1}^{n_k} x_{k;j} \right]^{w_k}\right) 
 \prod_{k=1}^m \prod_{j=1}^{n_k} dx_{k;j} \; , \label{f.def}
\ee 
where $\delta (x)$ represents Dirac's delta function that takes zero value everywhere except 
 at $x=0$ and that $\forall~a > 0$, $\int_{-a}^a \delta (x) dx = 1$. 

To proceed, let us make the following change of variables
\bea 
\tau &=& e^{-t} \nonumber \\
x_{k;j} &=& e^{-u_{k;j}} \nonumber
\eea
 and remember that if $y_0$ is the only root of $f$ ($f(y_0) = 0$) 
\[ 
 \delta (f(y)) = \frac{\delta(y-y_0)}{|f'(y_0)|} \; ,
\] 
we may then rewrite (\ref{f.def}) as 
\be
f(\tau)  = 
  \int_{0}^{\infty} \cdots \int_{0}^{\infty} e^t e^{-\sum_{k,j} u_{k;j}} 
 \delta \left( t - \sum_{k=1}^m w_k \left[ \sum_{j=1}^{n_k} u_{k;j} \right] \right) 
 \prod_{k=1}^m \prod_{j=1}^{n_k} du_{k;j} \; . \label{f.0}
\ee 
Note that the right hand side of (\ref{f.0}), except for the additional factor $e^t$, 
is the probability density function of a weighted, linear sum of exponential variables.  

By introducing the integral representation of the $\delta$ function
\[
\delta (t-c) = \frac{1}{2\pi} \int_{-\infty}^\infty dq \; e^{-iq(t-c)} \; ,
\]
we may re-express (\ref{f.0}) as 
\bea
f(\tau)  &=& \int_{-\infty}^\infty \frac{dq}{2\pi}\, e^{-it(q+i)}
\prod_{k=1}^m  \left[ \int_0^\infty  e^{-u} e^{i q w_k u} du  \right]^{n_k} \nonumber \\
&=& \int_{-\infty}^\infty \frac{dq}{2\pi}\, e^{-it(q+i)}
\prod_{k=1}^m \left[ \frac{1}{1-iq w_k}  \right]^{n_k} \nonumber \\
&=& \int_{-\infty}^\infty \frac{dq}{2\pi}\, e^{-it(q+i)}
 \prod_{l=1}^m \left( \frac{i}{w_l}\right)^{n_l}
\prod_{k=1}^m \left[\frac{1}{q + ir_k}  \right]^{n_k} \nonumber \\
&=& \left[ \prod_{l=1}^m  \, r_l^{n_l} \right] \left( i \right)^{\sum_{k=1}^m n_k}
\int_{-\infty}^\infty \frac{dq}{2\pi}\, e^{-it(q+i)}
\prod_{k=1}^m \left[\frac{1}{q + ir_k}  \right]^{n_k} \label{f.1} \\
& \equiv & \left[ \prod_{l=1}^m  \, r_l^{n_l} \right] \tilde f(\tau; n_1,n_2,\ldots,n_m) \;,
\label{tf.def}  
\eea 
where $r_k \equiv 1/w_k$ is introduced for the ease of analytical manipulation and $\tilde f$ is
 introduced for later convenience.  
Since all $w_k  > 0$, implying that all $r_k > 0$ and thus the poles of the integrand in 
(\ref{f.1}) lie completely at the lower half of the $q$-plane.  The integral of $q$ may be extended to
 enclose the lower half $q$-plane to result in 
\bea
f(\tau) &=&  e^t\, \left[ \prod_{l=1}^m   (i\, r_l)^{n_l} \right]\;  \left(\frac{-2\pi i}{2\pi} \right)
 \sum_{k=1}^m \, \frac{\left(\partial/ \partial q\right)^{n_k-1} }{(n_k-1)!}
 \left[\,  e^{-itq} \prod_{j=1, j\ne k}^m   \left(\frac{1}{q + ir_j}  \right)^{n_j}
 \right]_{q = -ir_k}  \nonumber \\
  &=& 
e^t\, \left[ \prod_{l=1}^m   (i\, r_l)^{n_l} \right]  \sum_{k=1}^m  \left\{ (-i) \!\!\!\!\!\!\!\!
\sum_{\substack{g_1,g_2,\ldots,g_m = 0\\ \sum g_i = n_k-1}} \!\!\!\! \frac{(-1)^{n_k-1} (it)^{g_k}}{g_k !\; e^{r_k\, t}} 
\prod_{j=1,j\ne k}^m \frac{ (n_j-1+g_j)!}{(n_j-1)!  g_j!} \left(\frac{-i}{r_j-r_k}  \right)^{n_j+g_j}
\right\}  \nonumber \\
  &=& 
e^t\, \left[ \prod_{l=1}^m  r_l^{n_l} \right]  \sum_{k=1}^m  \left\{ 
\sum_{\substack{g_1,g_2,\ldots,g_m = 0\\ \sum g_i = n_k-1}}\; \frac{ (t)^{g_k}}{g_k !} e^{-r_k\, t}
\prod_{j=1,j\ne k}^m \frac{ (n_j-1+g_j)!}{(n_j-1)!  g_j!} \frac{(-1)^{g_j}}{(r_j-r_k)^{n_j+g_j}}
\right\} 
 \; . \label{f.1.2}
\eea
Comparing eq.~(\ref{f.1.2}) with eq.~(\ref{f.0}), we see that the right hand side of (\ref{f.1.2}) 
 is composed of the product of the factor $e^t$ and the probability density function
  of a weighted sum of exponential variables. In fact, 
   the explicit expression of latter, in addition to the new derivation presented here, 
   was derived much earlier~\citep{Mathai_1983}  under different context and was rediscovered/rederived multiple    
    times~\citep{Harrison_1990,AM_1997,JK_2003} by different means. 
  Its connection to combining $P$-values, however, was never made explicit till now.

From (\ref{F.def}), we know that $F(\tau = 0 ) = 0$, implying that 
\bea
F(\tau) &=& \int_0^\tau f(\tau') d\tau' = \int_t^\infty f (e^{-t'})\; e^{-t'} dt' \nonumber \\
&=& \left[ \prod_{l=1}^m  r_l^{n_l} \right] \sum_{k=1}^m 
 \sum_{\substack{g_1,g_2,\ldots,g_m = 0\\ \sum g_i = n_k-1}}
\left(\prod_{j=1,j\ne k}^m \frac{ (n_j-1+g_j)!}{(n_j-1)!  g_j!} \frac{(-1)^{g_j}}{(r_j-r_k)^{n_j+g_j}} \right) 
 \int_t^\infty \frac{ (t')^{g_k}}{g_k !} e^{-r_k\, t'} dt' 
 \nonumber \\
&=&   \left[ \prod_{l=1}^m  r_l^{n_l} \right] \sum_{k=1}^m 
 \sum_{\substack{g_1,g_2,\ldots,g_m = 0\\ \sum g_i = n_k-1}}
\left(\prod_{j=1,j\ne k}^m \frac{ (n_j-1+g_j)!}{(n_j-1)!  g_j!} \frac{(-1)^{g_j}}{(r_j-r_k)^{n_j+g_j}} \right) 
\left( \sum_{l=0}^{g_k} \frac{t^{g_k-l}}{r_k^{l+1} (g_k-l)!} e^{-r_k\, t}\right)
 \nonumber \\
&=&  \sum_{k=1}^m 
 \sum_{\substack{g_1,g_2,\ldots,g_m = 0\\ \sum g_i = n_k-1}}
\left(\prod_{j=1,j\ne k}^m \frac{ (n_j-1+g_j)!}{(n_j-1)!  g_j!} \frac{(-r_k)^{g_j} r_j^{n_j}}{(r_j-r_k)^{n_j+g_j}} \right) 
\left( \sum_{l=0}^{g_k} \frac{(r_k\, t)^{g_k-l}}{(g_k-l)!} e^{-r_k\, t}\right)
 \nonumber \\
&=&  \sum_{k=1}^m \sum_{g_k=0}^{n_k-1} 
 \sum_{\substack{g_{i\ne k} = 0\\ \sum_{i} g_i = n_k-1}}^{n_k-1-g_k}
\left(\prod_{j=1,j\ne k}^m \frac{ (n_j-1+g_j)!}{(n_j-1)!  g_j!} \frac{(-r_k)^{g_j} r_j^{n_j}}{(r_j-r_k)^{n_j+g_j}} \right) 
 H(r_k\, t,\, g_k) \;,  \label{F.unify.0}
\eea
where the function $H$ is defined as 
\be \label{H.def}
H(x,n) \, \equiv \;  e^{-x}\, \sum_{k=0}^{n} \frac{x^k}{k!}  \; .
\ee
Eq.~(\ref{F.unify.0}) represents the most general formula that
interpolates the scenarios considered by Fisher and Good.

Let us take the limiting case from (\ref{F.unify.0}). For Fisher's formula, one weights every $P$-value equally,
 and thus correspond to $m=1$ and $n_1 = L$. The constraint in the sum of (\ref{F.unify.0}) forces 
$g_1 = n_1-1 = L-1$. Consequently, we have (by calling $r_1$ by $r$ for simplicity)
\be 
\Prob(Q_F\le \tF) = H(rt,\, L-1)
 =   e^{-rt} \, \sum_{l=0}^{L-1} \frac{(rt)^{l}}{l!}
  \label{Unify.Fisher.0}
\ee 
Notice that regardless whatever the weight $w$ one assigns to all the $P$-values, the final answer is independent of the 
 weight. This is because $t = -\ln \tau = -w \ln \tF =  (-\ln \tF)/r $ and therefore $rt = \ln (1/\tF)$.
 This results in     
\be  \label{F.unify.Fisher}
\Prob(Q_F\le \tF)  =  \tF \sum_{l=0}^{L-1} \frac{\left[\; \ln (1/\tF)\, \right]^{\, l} }{l!} \;, 
\ee
exactly what one anticipates from (\ref{Comb.Fisher}). 
To obtain the results of Good, one simply makes $m=L$ and $n_k = 1$ $\forall~k$, implying all $g_i = 0$. In this case, (\ref{F.unify.0}) becomes  (with $r_l = 1/w_l$, $e^{-t} = \tG$ and $H(x,0)=1$ $\forall~x$)  
 \be  \label{F.unify.Good}
\Prob (Q_G \le \tG) = \sum_{l=1}^L  \left( \prod_{k \ne l} \frac{r_k}{r_k-r_l} \right) \tG^{1/w_l} 
\;  = \;  \sum_{l=1}^L \Lambda_l \; \tG^{1/w_l} \; , 
 \ee
reproducing exactly (\ref{Comb.Good})

One may also re-express eq.~(\ref{F.unify.0}) in a slightly different form
\bea 
F(\tau) &=& \left[ \prod_{l=1}^m  r_l^{n_l} \right] \sum_{k=1}^m \sum_{g_k=0}^{n_k-1} 
\frac{1}{r_k^{g_k+1}}  H(r_k\, t,\, g_k) \times \nonumber \\
&& \hspace*{30pt} \times \sum_{\substack{g_{i\ne k} = 0\\ \sum_{i} g_i = n_k-1}}^{n_k-1-g_k}
\left(\prod_{j=1,j\ne k}^m \frac{ (n_j-1+g_j)!}{(n_j-1)!  g_j!} \frac{(-1)^{g_j}}{(r_j-r_k)^{n_j+g_j}} \right) \nonumber \\
&\equiv & \left[ \prod_{l=1}^m  r_l^{n_l} \right] \tilde F (\tau ;  n_1, n_2, \ldots, n_m) \; .
\label{F.unify.1}
\eea
Note that in the expression (\ref{F.unify.1}), we have isolated an overall multiplying factor 
 and keeps explicit the $n_{1\le k \le m}$ dependence for later convenience.  
%is the main result of this paper. It represents the most general formula that
% interpolates the scenarios considered by Fisher and Good. 

\section{Cases of nearly degenerate weights} 
In our derivation of (\ref{F.unify.Fisher}), it is explicitly shown in section~\ref{sec:derivation} 
 that the final $P$-value obtained is 
 independent of the weight $w$ that was used to assign to all the individual $P$-values, $p_1,p_2,\ldots,p_L$.  
 It is thus natural to ask, if one starts by weighting each $P$-value differently, upon making 
  the weights closer to one another, will one recover Fisher's formula (\ref{Comb.Fisher})
   from Good's formula (\ref{Comb.Good})?  By continuity, the answer is expected to be affirmative. 
    More generally, one would like to have a formal protocol to
    compute the combined $P$-value when the weights may be categorized into several subsets, 
    within each subset the weights are almost degenerate. 
   
 In this section, we first illustrate the transition from Good's formula to
  Fisher's formula by combining two almost degenerate $P$-values. We will then provide a general
   protocol to explicitly, when there exist nearly degenerate weights, deal with the possible
    numerical instability that was first cautioned by~\citet{Good_1955} and subsequently by 
   many authors~\citep{Bhoj_1992,OS_2001,Hou_2005}.
    
Let us consider combining $p_1$ and $p_2$ with weights $w_1$ and $w_2$ using Good's formula. 
One has
\be \label{Good.2p.0}
\Prob (Q_G \le \tG) = \frac{1}{w_1-w_2} \left[ w_1\, p_1 p_2^{\frac{w_2}{w_1}} 
- w_2\, p_1^{\frac{w_1}{w_2}} p_2  \right] \; .
\ee  
Without loss of generality, one assumes $w_1 > w_2$ and hence writes $w_1/w_2 = 1+ \epsilon$ with
 $\epsilon > 0$.  We are interested in the case when the weights get close to each other, or when
  $\epsilon \to 0$. We now rewrite eq.~(\ref{Good.2p.0}) as 
\be \label{Good.2p.1}
\Prob (Q_G \le \tG) = \frac{w_2}{w_1-w_2} \left[ \frac{w_1}{w_2}\, p_1 p_2^{\frac{w_2}{w_1}} 
- \, p_1^{\frac{w_1}{w_2}} p_2  \right] = \frac{1}{\epsilon}
  \left[ (1+\epsilon) \, p_1 p_2^{\frac{1}{1+\epsilon}} 
- \, p_1^{1+\epsilon} p_2  \right] 
 \; .
\ee   
In the limit of small $\epsilon$, we may rewrite (\ref{Good.2p.1})as 
\bea
\Prob (Q_G \le \tG) &=&  \frac{p_1p_2}{\epsilon}
  \left[ (1+\epsilon) \, p_2^{-\frac{\epsilon}{1+\epsilon}} 
- \, p_1^{\epsilon}   \right]  = \frac{p_1p_2}{\epsilon}
  \left[ (1+\epsilon) \, e^{-\frac{\epsilon}{1+\epsilon} \ln p_2} 
- \, e^{\epsilon \ln p_1}   \right] \nonumber \\
&=& \frac{p_1p_2}{\epsilon}  \left[ \epsilon -\epsilon (\ln p_2 + \ln p_1 )  + {\cal O}(\epsilon^2)
% + \frac{\epsilon^2}{2}\left( \frac{1}{1+\epsilon} \ln^2 p_2 - \ln^2 p_1\right)
\right] \nonumber \\
&=& p_1p_2  \left[ 1 -\ln (p_1 p_2) + {\cal O}(\epsilon)
 %+ \frac{\epsilon}{2}\left( \frac{1}{1+\epsilon} \ln^2 p_2 - \ln^2 p_1\right) 
\right] \nonumber \\
&\xrightarrow[\epsilon \to 0]{} & p_1p_2  \left[ 1 -\ln (p_1 p_2) \right] \; = \; \Prob (Q_F \le \tF) \; .
 \label{Good.2p.to.Fisher}
\eea
Note that when the small weight difference $w_1-w_2$ is near the machine precision of a digital computer, 
 using formula (\ref{Comb.Good}) directly will inevitably introduce numerical instability caused by 
 rounding errors.
% The series expansion shown in (\ref{Good.2p.to.Fisher}) above then provides a good way for 
%  accurately computing the final $P$-value.  
  
To construct a general protocol to deal with nearly degenerate weights, one first observes 
 from eqs.~(\ref{f.1}-\ref{F.unify.1}) that it is the inverse of weights $r_k \equiv 1/w_k$ that permeates
 the derivation of the unified $P$-value. The closeness between weights is thus naturally defined by closeness  
 in the inverse weights. As shown in eqs.~(\ref{Good_tau}) and (\ref{Comb.Good}), the combined $P$-value 
  by Good's formula is independent of the absolute size of the weights but only on the relative weights.
   Making the observation that $r_k\, t$ in eq.~(\ref{F.unify.0}) only depend on the ratios $r_{j\ne k}/r_k$, one also
    sees explicitly that the most general combined $P$-values (see (\ref{F.unify.0})) only depend on 
     the relative weights as well. We are thus free to choose any scale we wish.
 For simplicity, we normalize the inverse weight associated with each method by demanding 
  the sum of inverse weights equal the total number of methods 
\be  \label{iw.norm}
\sum_{j=1}^M r_j = \sum_{j=1}^M 1 = M\; ,
\ee
 where $1/r_j$ represents the weight associated method $j$ and $M$ represents the total
  number of $P$-values (or methods) to be combined. For the GC described in section~\ref{sec:derivation}, 
  $M = \sum_{k=1}^m n_k$.  
  This normalization choice makes the average inverse weight of participating methods be $1$.  

The next step is to determine, for a given list of inverse weights and 
 the radius of clustering, the number of clusters needed.   
This task may be achieved in a hierarchical manner. After normalizing the inverse weights $r_k$
 using  eq.~(\ref{iw.norm}), one may sort the inverse weights in either ascending or descending order.
 For a given radius $\eta > 0$, one starts to seek the pair of inverse weights that are closest but not identical, 
 and check if it is smaller than the radius $\eta$. If yes, one will merge that pair of inverse weights by using their  
  average, weighted by number of occurrences, as the new center and continue the process till every inverse weights in the
   list is separated by a distance farther than $\eta$. 
   We use an example of $M=8$ to illustrate the idea. Let the normalized inverse weights $\{r_j\}_{j=1}^8$ be
 \[
 0.50,~0.70(2),~0.71,~0.74,~1.03~,1.80~,1.82
 \]  
 where the number $2$ inside the pair of parentheses after $0.70$ simply indicates that there
  are two identical inverse weights $0.70$ to start with. Assume that one chooses the radius of
   cluster $\eta$ to be $0.005$, since every adjacent inverse weights are separated by more than
    $0.005$, no further clustering procedures is needed and one ends up having seven
     effective clusters: one cluster with two identical inverse weights $0.70$, and the rest of six clusters
   are all singletons. This corresponds to $m=7$, $n_1 =1$, $n_2 = 2$, $n_3 = n_4 =\cdots = n_7 =1 $.

  Suppose one chooses the clustering radius
  $\eta$ to be $0.05$. In the first step, we identify that $0.70$ and $0.71$ are the closet pair of 
  inverse weights. The weighted average between them is
 \[
 \frac{2 \cdot 0.70 + 0.71}{3} = \frac{2.11}{3}= 0.70\bar3 \; .
 \]    
The list of inverse weights then appears as
\[
0.50,~0.70\bar 3(3),~0.74,~1.03~,1.80~,1.82 \; .
\]
The closest pair of inverse weights is now between $1.80$ and $1.82$, and upon merging them
 we have the list now appears as
\[
0.50,~0.70\bar 3(3),~0.74,~1.03~,1.81(2) \; .
\]
Next pair of closest inverse weights is then $0.70\bar 3$ and $0.74$.
 The weighted average leads to $(2.11+0.74)/4 = 0.7125$. 
After this step, the list of inverse weights appears as
\[
0.50,~0.7125(4),~1.03~,1.81(2) \; ,
\]    
indicating that we have $m=4$ ( four clusters), with number of members being $n_1 = 1$, $n_2 = 4$, $n_3 = 1$ and $n_4 = 2$. 
The centers of the four clusters are specified by the average inverse weights: $0.50,~0.7125,~1.03~,1.81$.
 The distance between any two of the average inverse weights is now larger than $0.05$. 
 
 This is a good place for us to introduce some notation. We shall denote by $r_{k} + \eta_{k;j}$ the $j$th 
 inverse weights of cluster $k$, whose averaged inverse weight is $r_k$. With this definition, for the example above, 
 we have $\eta_{1;1} = 0$, 
 $\eta_{2;1} = \eta_{2;2} = -0.0125$, $\eta_{2;3} = -0.0025$, $\eta_{2;4} = 0.0275$, $\eta_{3;1} = 0$, 
 $\eta_{4;1} = -0.01$, and $\eta_{4;2} = 0.01$.

Using the hierarchical protocol mentioned above, the number of clusters $m$ and the numbers of members $n_k$
 associated with cluster $k$ are all obtained along with $\{ \eta_{k;j} \}$. Following the derivation in 
 section~\ref{sec:derivation}, we obtain a probability density function very similar to  (\ref{f.1})  
\be   
f(\tau)  =  \left[ \prod_{l=1}^m  \prod_{j=1}^{n_l} ( r_l + \eta_{l;j}) \right] \left( i \right)^{\sum_{k=1}^m n_k}
\int_{-\infty}^\infty \frac{dq}{2\pi}\, e^{-it(q+i)}
\prod_{k=1}^m \left[ \prod_{j=1}^{n_k} \frac{1}{q + i(r_k + \eta_{k;j})}  \right]  \; .\label{f.2}
% \nonumber \\
%&=& \left[ \prod_{l=1}^m \left( ir_l \right)^{n_l} e^{\sum_{j=1}^{n_l} \ln (1+ \frac{\eta_{l;j}}{r_l})} \right]
%\int_{-\infty}^\infty \frac{dq}{2\pi}\, e^{-it(q+i)}
%\prod_{k=1}^m \left\{ \left(\frac{1}{q + ir_k } \right)^{n_k}  \right. \nonumber \\
%&& \hspace*{180pt}\prod_{j=1}^{n_k}  \left. \left(1+ \frac{i\eta_{k;j}}{q + ir_k } \right) \right\} 
\ee
%with $r_k$ representing the average inverse weights of the $k$th cluster and
%$\eta_{k;j}$ representing the difference between the $j$th inverse weight within cluster $k$ and 
%$r_k$.  
    
From section~\ref{sec:derivation}, we see that the ill-conditioned situations emerge 
 when some weights are nearly degenerate and the source of difference in inverse weights comes from
  obtaining $\tilde F(\tau; n_1,n_2,\ldots,n_m)$ in (\ref{F.unify.1}) 
from $\tilde f (\tau ; n_1,n_2,\ldots,n_m)$ in (\ref{tf.def}).  Therefore, one may leave the
 prefactor $\left[ \prod_{l=1}^m  \prod_{j=1}^{n_l} ( r_l + \eta_{l;j}) \right]$ untouched and focus on
 the rest of the right hand side of eq.~(\ref{f.2}). To proceed, we write 
\bea 
\frac{1}{q + i(r_k + \eta_{k;j})} &=& \frac{1}{q + ir_k } \left(1 + \frac{i\,\eta_{k;j}}{q + ir_k} \right)^{-1} 
= \frac{1}{q + ir_k }\, e^{-\ln \left(1 +  \frac{i\,\eta_{k;j}}{q + ir_k} \right) } \nonumber \\
&=& \frac{1}{q + ir_k } \, \exp \left[\sum_{g=1}^\infty \frac{1}{g}
\left(\frac{-i\,\eta_{k;j}}{q + ir_k}\right)^g \right]  \; .
\nonumber
\eea 
Consequently, we may write
\be \label{Ykg.intro} 
\prod_{j=1}^{n_k} \frac{1}{q + i(r_k + \eta_{k;j})} =
\frac{1}{(q + ir_k)^{n_k}} \exp \left[\sum_{g=1}^\infty \frac{ Y_{k;g} \, (i)^g}
{\left(q + ir_k\right)^g} \right] 
\ee     
where
\be \label{Ykg.def} 
Y_{k;g} \equiv \sum_{j=1}^{n_k} \frac{(-\eta_{k;j})^g}{g} \; .
\ee
The product in eq.~(\ref{f.2}) may now be formally written as 
\be \label{prod.0}
\prod_{k=1}^m \left[ \prod_{j=1}^{n_k} \frac{1}{q + i(r_k + \eta_{k;j})}  \right]
 =  \left[ \prod_{k=1}^m \frac{1}{(q + ir_k)^{n_k}}  \right] 
 \exp \left[ \sum_{g=1}^\infty (i)^g \sum_{k=1}^m \frac{ Y_{k;g} }
{\left(q + ir_k\right)^g} \right] \; .
\ee

We now note a simplification by choosing $r_k$ to be the average inverse weight of the $k$th cluster.
 In this case, we have $\sum_{j=1}^{n_k} \eta_{k;j} = 0$ $\forall~k$. That is, $Y_{k;1} = 0$ always.
 This allows us to write eq.~(\ref{prod.0}) as 
\be \label{prod.1} 
\prod_{k=1}^m \left[ \prod_{j=1}^{n_k} \frac{1}{q + i(r_k + \eta_{k;j})}  \right]
 =  \left[ \prod_{k=1}^m \frac{1}{(q + ir_k)^{n_k}}  \right] 
 \exp \left[ \sum_{g=2}^\infty (i)^g \sum_{k=1}^m \frac{ Y_{k;g} }
{\left(q + ir_k\right)^g} \right] \; .
\ee 
        
The key idea here is to Taylor expand the exponential and collect terms of 
 equal number of $1/(q+ir)$. Evidently, the first correction term starts with $1/(q+ir)^2$. 
 Furthermore, before the $1/(q+ir)^4$ order, there is no mixing between different clusters.
Below, we rewrite eq.~(\ref{f.2}) to include the first few orders of correction terms
\bea
&& \hspace*{-40pt}\frac{f(\tau)}{\prod_{l=1}^m  \prod_{j=1}^{n_l} ( r_l + \eta_{l;j})}  
=  %\left[ \prod_{l=1}^m  \prod_{j=1}^{n_l} ( r_l + \eta_{l;j}) \right] 
\left( i \right)^{\sum_{k=1}^m n_k}
\int_{-\infty}^\infty \frac{dq}{2\pi}\, e^{-it(q+i)} 
\frac{ \exp \left[ \sum_{g=2}^\infty (i)^g \sum_{k=1}^m \frac{ Y_{k;g} }
{\left(q + ir_k\right)^g} \right]}{\prod_{k=1}^m (q + ir_k)^{n_k}} \nonumber \\
&& = \tilde f(\tau; \{ n_l \}_{l=1}^m) + \sum_{k=1}^m Y_{k;2}\; \tilde f(\tau; \{ n_{l\ne k}, n_k+2 \}) \nonumber \\
&& +  \sum_{k=1}^m Y_{k;3}\; \tilde f(\tau; \{ n_{l\ne k}, n_k+3 \}) 
+  \sum_{k=1}^m \left(Y_{k;4} + \frac{(Y_{k;2})^2}{2!} \right)\; \tilde f(\tau; \{ n_{l\ne k}, n_k+4 \}) \nonumber \\
&& +  \frac{1}{2!}\sum_{\substack {k, k' =1 \\ k\ne k'}}^m Y_{k;2} Y_{k';2}\; \tilde f(\tau; \{ n_{l\ne k,k'}, n_k+2,n_{k'}+2 \})
 + {\cal O}(\eta^5) \; .
\eea   
This immediately leads to 
\bea 
&& \hspace*{-40pt}\frac{F(\tau)}{\prod_{l=1}^m  \prod_{j=1}^{n_l} ( r_l + \eta_{l;j})} 
  = \tilde F(\tau; \{ n_l \}_{l=1}^m) + \sum_{k=1}^m Y_{k;2}\; \tilde F(\tau; \{ n_{l\ne k}, n_k+2 \}) \nonumber \\
&& +  \sum_{k=1}^m Y_{k;3}\; \tilde F(\tau; \{ n_{l\ne k}, n_k+3 \}) 
+  \sum_{k=1}^m \left(Y_{k;4} + \frac{(Y_{k;2})^2}{2!} \right)\; \tilde F(\tau; \{ n_{l\ne k}, n_k+4 \}) \nonumber \\
&& +  \frac{1}{2!}\sum_{\substack {k, k' =1 \\ k\ne k'}}^m Y_{k;2} Y_{k';2}\; \tilde F(\tau; \{ n_{l\ne k,k'}, n_k+2,n_{k'}+2 \})
 + {\cal O}(\eta^5) \; . \label{F.unify.expansion}
\eea 

Note that when the cluster radius $\eta$ is chosen to be zero, the only clusters are from 
 sets of {\it identical} weights, and all $\eta_{k;j}$ must be zero. In this case,
  only the first term on the right hand side of (\ref{F.unify.expansion}) exists and 
  the result derived in section~\ref{sec:derivation} is recovered exactly.  
Since all $\tilde F$ are finite positive quantities, the errors resulting from truncating 
 the expression in eq.~(\ref{F.unify.expansion}) at certain
 order of $\eta$ can be easily bounded. Therefore, any desired precision may be obtained via 
 including the corresponding number of  higher order terms. As the main result of the
   current paper, our expansion provides a systematic, numerically stable method 
   to achieve desired accuracy in computing combined $P$-values.

\section{Examples}

{\bf Example~(a)}: This example, assuming $m=4$, demonstrates how to 
compute the $\tilde F(\tau;\{n_l \})$ function present in eq.~(\ref{F.unify.1}).
 Let $r_k$ be the inverse weights associated with cluster $k$. 
 When combining multiple $P$-values with a prescribed clustering radius on the inverse weights, see (\ref{iw.norm})
 and the previously described clustering procedure, the $r_k$s and the deviations 
  $\eta_{k;j}$ are obtained once and for all. The $\eta_{k;j}$, through eq.~(\ref{Ykg.def}),  
  constitute the key expansion parameters $Y_{k;g}$ that yield, upon multiplying 
   by $\tilde F(\tau;{n_l})$ with different $\{ n_l \}$, the higher order terms in our key result (\ref{F.unify.expansion}).   
 Note that in eq.~(\ref{F.unify.expansion}), in the zeroth order term, the argument $n_l$ of $\tilde F$  
  represents the number of members associated with cluster $l$. However, for higher order correction terms,  
 the $n_l$s entering $\tilde F$ no longer carry the same meaning.  Therefore, in the example shown
  here, one does not assume that $n_j$ is the number of methods associated with cluster $k$.  
 We now illustrate how to open up the sum in eq.~(\ref{F.unify.1}). 
  The constraint $\sum_i g_i = n_k -1$ implies that one only has $m-1$ independent $g_i$s.   
  Once $m-1$ $g_i$s are specified, the remaining one is also determined.  
  To simplify the exposition, let us introduce the following notation
  \[
 \alpha(g_j;j,k) \equiv  \frac{ (n_j-1+g_j)!}{(n_j-1)! g_j!} \frac{(-1)^{g_j}}{(r_j-r_k)^{n_j+g_j}} \; .
  \]

One then expand the sum in (\ref{F.unify.1}) as 
\bea
 \tilde F(\tau) = \sum_{g_1=0}^{n_1-1}  \frac{ H(r_1\, t,\, g_1) }{r_1^{g_1+1}} 
 \sum_{g_2=0,}^{n_1-1-g_1} \alpha(g_2;2,1) \; \sum_{g_3=0}^{n_1-1-g_1-g_2} \alpha(g_3;3,1)\; \; \alpha(g_4;4,1) \nonumber \\ 
+ \sum_{g2=0}^{n_2-1} \frac{ H(r_2\, t,\, g_2)}{r_2^{g_2+1}}  \sum_{g_1=0}^{n_2-1-g_2} \alpha(g_1;1,2) 
\;\sum_{g_3=0}^{n_2-1-g_2-g_1}  \alpha(g_3;3,2)\; \; \alpha(g_4;4,2)  \nonumber \\ 
+ \sum_{g3=0}^{n_3-1}  \frac{H(r_3\, t,\, g_3)}{r_3^{g_3+1}}   \sum_{g_1=0}^{n_3-1-g_3} \alpha(g_1;1,3) 
\;\sum_{g_2=0}^{n_3-1-g_3-g_1}  \alpha(g_2;2,3)\; \; \alpha(g_4;4,3)  \nonumber \\ 
+ \sum_{g4=0}^{n_4-1} \frac{ H(r_4\, t,\, g_4)}{r_4^{g_4+1}}  \sum_{g_1=0}^{n_4-1-g_4} \alpha(g_1;1,4) 
\;\sum_{g_2=0}^{n_4-1-g_4-g_1}  \alpha(g_2;2,4)\; \; \alpha(g_3;3,4)  \; .
\label{example:a}
\eea 

%which constitute the contribution to the zero order term of equation~\ref{F.unify.expansion}. 
\vspace*{24pt}

\noindent {\bf Example~(b)}: This example illustrates the possibility of numerical instability 
associated with eqs.~(\ref{Comb.Good}) and (\ref{F.unify.1}) when they are used to combine 
{\it P}-values with nearly equal weights. We also show how such instabilities  are resolved 
by using eq.~(\ref{F.unify.expansion}).
Consider the case of combining five $P$-values, $\{0.008000257, 0.008579261,  0.0008911761,  0.006967988,  0.004973110\}$, 
 weighted respectively by $\{0.54531152, 0.54532057, 0.54531221, 0.54531399,  0.54531776\}$. 
 % Upon normalizations, the corresponding inverse weights are given by $\{ XXXXX \}$, leading to 
 Using eq.~(\ref{Good_tau}), one obtains $\tau_G =  4.30656196 \times 10^{-7}$. 
 The combined $P$-value is then obtained as the probability of  attaining
  a random variable $Q_G$, defined in eq.~(\ref{Good_Q}), such that is less than or equal to $\tau_G$. 

Combining $P$-values using eq.~(\ref{Comb.Good}) gives

\bea
&\Prob (Q_G \le \tau_G)& = 1923475672.53812003 +  134195847.49348195 \nonumber \\
& &  - 3271698577.16100168 + 1726093852.57087326 - 512066795.44147670 \nonumber  \\
& &=   -0.00000322  \; . \nonumber
\eea
When one uses equation~(\ref{F.unify.1}), $\tau$ takes the value of $\tau_G$ and 
 the random variable $Q$ is simply $Q_G$, and the combined $P$-value becomes
\bea
&\Prob (Q \le \tau)& = 170090507.09336647 + 21761086.68190728   \nonumber \\
& & - 972903041.25101399 + 941269625.31004059 - 512066795.44252247  \nonumber \\
& &=   -0.00000006   \; .\nonumber
\eea
Apparently,  probability can't be negative and the negative values shown above 
illustrate how eqs.~(\ref{Comb.Good}) and (\ref{F.unify.1}) may suffer from numerical instability
 when the weights are nearly degenerate.  
 This numerical instability is removed by applying equation~(\ref{F.unify.expansion}) which combines 
 weighted $P$-values using a controlled expansion and yields, for this example, 
\bea
&\Prob (Q \le \tau )& =  5.379093 \times 10^{-8} +   1.407305\times 10^{-16}  \nonumber \\
& &  -1.066323 \times 10^{-21} +  1.634917 \times 10^{-25}    + {\cal O}(10^{-29}) \nonumber \\
& &=   5.37909 \times 10^{-8}  \nonumber
\eea

\vspace*{24pt}

\noindent {\bf Example~(c)}: 
 One natural question to ask is that how well does eq.~(\ref{F.unify.expansion}) work
  when one chooses a larger clustering radius and group weights that are clearly distinguishable into
   a few clusters? To consider this case, let us use the five $P$-values from example (b) above
   but with weights chosen differently. Let us assume that the inverse weights ($r_k \equiv 1/w_k $) 
   associated with these five $P$-values are  $\{ 0.6, 0.65, 1.2, 1.25,1.3\}$. For this case, 
  $\tau = \tau_G =1.935663\times 10^{-13}$.  Combining $P$-value using    
    formulas (\ref{Comb.Good}) yields 
 \bea 
\Prob (Q \le \tau) &=& \;\; 2.187324 \times 10^{-6} - 5.946040 \times 10^{-7} + 2.131226 \times 10^{-13}\nonumber \\
&&  - 8.011644 \times 10^{-14} + 7.639290 \times 10^{-15}\nonumber \\
&=& 1.59272 \times 10^{-6} \; ,  \nonumber 
 \eea    
while combining $P$-values using (\ref{F.unify.1}) yields identical results 
 \bea 
\Prob (Q \le \tau) &=& \;\; 1.725699 \times 10^{-6} -3.049251 \times 10^{-7} + 1.311524 \times 10^{-13}\nonumber \\
&&  -6.162803 \times 10^{-14} + 7.639290 \times 10^{-15}\nonumber \\
&=& 1.59272 \times 10^{-6} \; .  \nonumber 
 \eea  

When one uses $\eta = 0.1$ as the clustering radius, one obtains two clusters: one with average inverse weight
 $0.625$ and the other with average inverse weight $1.25$. 
 If one then uses eq.~(\ref{F.unify.expansion}) to combine $P$-values, 
  one attains the following results     
\bea 
\Prob (Q \le \tau)  &=& \;\; 1.472453 \times 10^{-6} + 1.171521 \times 10^{-7} + 0 \nonumber \\
&&  + 2.584710 \times 10^{-9} + 4.889899 \times 10^{-10} + {\cal O}(10^{-12})\nonumber \\
&=& 1.59268 \times 10^{-6} \; ,  \label{Example.exp} 
\eea 
which contains no sign alternation and 
agrees well with the results from both (\ref{Comb.Good}) and (\ref{F.unify.1}).     
This illustrates the robustness of eq.~(\ref{F.unify.expansion}) in combining $P$-values.    
Note that the third term on the right hand side of (\ref{Example.exp}) is zero. 
 This is because the multiplying factor $Y_{k;3}$ is zero for both clusters. 
 In general, $Y_{k;3}$ measures the skewness of inverse weights associated with cluster $k$ and
  for our case here both cluster of inverse weights are perfectly symmetrical with respect to
   their centers, leading to zero skewness. If the inverse weights of 
   cluster $k$ distribute perfectly symmetrically with respect to its center, 
   it is evident from eq.~(\ref{Ykg.def}) that $Y_{k;g} = 0$ for odd $g$.

Evidently if one chooses a large clustering radius $\eta$ and then uses eq.~(\ref{F.unify.expansion})
 to combine $P$-values, many higher order terms in the expansion will be required 
 to achieve high accuracy in the final combined $P$-value.

\section{Future directions} 
Although the expression (\ref{F.unify.0}) provides access to exact statistics for a 
 broader domain of problems and our expansion formula (\ref{F.unify.expansion}) 
 provides accurate and stable statistics even when nearly degenerate weights are present, 
 there remain a few unanswered questions that should be addressed by the community in the near future. 
 For example, even though we can accommodate any reasonable $P$-value weighting, thanks to (\ref{F.unify.expansion}), 
 the more difficult question is how does one choose the right set of weights when combining statistical significance~\citep{ZJ_1959,KP_1978,HO_1985,PF_1989,Forrest_2001}.
  The weights chosen reflects how much does one wish to trust various obtained $P$-values. 
  Ideally, a fully systematic method should also provide a metric for choosing appropriate weights.
  How to obtain the best set of weights remains an open problem and definitely deserves further 
  investigations. 
  
Another limitation of (\ref{F.unify.0}) and (\ref{F.unify.expansion}), 
and consequently of Fisher's and Good's formulas, is that one must assume the $P$-values to 
be combined as independent. In real applications,
it is foreseeable that $P$-values reported by various methods may exhibit non-negligible 
correlations. How to obtain the correlation~\citep{WL_1984,PGT_1987,James_1991} 
and properly take into account the existence of $P$-value correlations~\citep{Brown_1975,KM_2002,Hou_2005}
is also a challenging problem that we hope to address in the near future.

\section*{Acknowledgement} 
This work was supported by the Intramural Research Program of the National Library of Medicine at National Institutes of Health.
% This study utilized the high-performance computational capabilities of the Biowulf Linux cluster at the National Institutes of Health, Bethesda, MD. (http://biowulf.nih.gov).


\begin{thebibliography}{}

\bibitem[\protect\citeauthoryear{Alves, Wu, Wang, Shen, and Yu}{Alves
  et~al.}{2008}]{AWWSetal_2008}
Alves, G., W.~W. Wu, G.~Wang, R.~F. Shen, and Y.~K. Yu (2008, Aug).
\newblock {{E}nhancing peptide identification confidence by combining search
  methods}.
\newblock {\em J. Proteome Res.\/}~{\em 7}, 3102--3113.

\bibitem[\protect\citeauthoryear{Amari and Mirsa}{Amari and
  Mirsa}{1997}]{AM_1997}
Amari, S. and R.~B. Mirsa (1997).
\newblock {Closed-form expression for distribution of the sum of independent
  exponential random variables}.
\newblock {\em IEEE Trans. Reliability\/}~{\em 46\/}(4), 519--522.

\bibitem[\protect\citeauthoryear{Bahrucha-Reid}{Bahrucha-Reid}{1960}]{Bharucha%
-Reid_1960}
Bahrucha-Reid, A. (1960).
\newblock {\em Elements of the Theory of Markov Processes and their
  Applications.}
\newblock McGraw-Hill.

\bibitem[\protect\citeauthoryear{Bailey and Gribskov}{Bailey and
  Gribskov}{1998}]{BG_1998}
Bailey, T.~L. and M.~Gribskov (1998).
\newblock {{C}ombining evidence using p-values: application to sequence
  homology searches}.
\newblock {\em Bioinformatics\/}~{\em 14}, 48--54.

\bibitem[\protect\citeauthoryear{Bateman}{Bateman}{1910}]{Bateman_1910}
Bateman, H. (1910).
\newblock A solution of a system of differential equations occurring in the
  theory of radiactive transformation.
\newblock {\em Pro. Cambridge Philosophical Soc.\/}~{\em 15}, 423--427.

\bibitem[\protect\citeauthoryear{Bhoj}{Bhoj}{1992}]{Bhoj_1992}
Bhoj, D.~S. (1992).
\newblock On the distribution of the weighted combination of independent
  probabilities.
\newblock {\em Statistics \& Probability Letters\/}~{\em 15\/}(1), 37 -- 40.

\bibitem[\protect\citeauthoryear{Brown}{Brown}{1975}]{Brown_1975}
Brown, M.~B. (1975).
\newblock A method for combining non-independent, one-sided tests of
  significance.
\newblock {\em Biometrics\/}~{\em 31\/}(4), 987--992.

\bibitem[\protect\citeauthoryear{Cox}{Cox}{1962}]{Cox_1962}
Cox, D. (1962).
\newblock {\em Renewal Theory}.
\newblock Methuen \& Co.

\bibitem[\protect\citeauthoryear{Feller}{Feller}{1966}]{Feller_1966}
Feller, W. (1966).
\newblock {\em An Introduction to Probability Theory and Its Applications},
  Volume~2.
\newblock John Wiley \& Sons.

\bibitem[\protect\citeauthoryear{Fisher}{Fisher}{1932}]{Fisher_1932}
Fisher, R.~A. (1932).
\newblock {\em Statistical Methods for Research Workers, vol. II}.
\newblock Edinburgh: Oliver and Boyd.

\bibitem[\protect\citeauthoryear{Forrest}{Forrest}{2001}]{Forrest_2001}
Forrest, W.~F. (2001).
\newblock Weighting improves the "new {H}aseman-{E}lston" method.
\newblock {\em Hum Hered\/}~{\em 52\/}(1), 47--54.

\bibitem[\protect\citeauthoryear{Gabler}{Gabler}{1987}]{Gabler_1987}
Gabler, S. (1987).
\newblock A quick and easy approximation to the distribution of a sum of
  weighted chi-square variables.
\newblock {\em Statistische Hefte\/}~{\em 28}, 317--25.

\bibitem[\protect\citeauthoryear{Good}{Good}{1955}]{Good_1955}
Good, I.~J. (1955).
\newblock On the weighted combination of significance tests.
\newblock {\em Journal of the Royal Statistical Society. Series B
  (Methodological)\/}~{\em 17\/}(2), 264--265.

\bibitem[\protect\citeauthoryear{Harrison}{Harrison}{1990}]{Harrison_1990}
Harrison, P.~G. (1990).
\newblock Laplace transform inversion and passage-time distribution in markov
  processes.
\newblock {\em J. Appl. Prob.\/}~{\em 27\/}(1), 74--87.

\bibitem[\protect\citeauthoryear{Hedges}{Hedges}{1985}]{HO_1985}
Hedges, L.~V., O.~I. (1985).
\newblock {\em Statistical Methods for Meta-Analysis}.
\newblock Academic Press.

\bibitem[\protect\citeauthoryear{Hess and Iyer}{Hess and Iyer}{2007}]{HI_2007}
Hess, A. and H.~Iyer (2007).
\newblock Fisher's combined p-value for detecting differentially expressed
  genes using {A}ffymetrix expression arrays.
\newblock {\em BMC Genomics\/}~{\em 8\/}(NIL), 96.

\bibitem[\protect\citeauthoryear{Hou}{Hou}{2005}]{Hou_2005}
Hou, C.-D. (2005).
\newblock A simple approximation for the distribution of the weighted
  combination of non-independent or independent probabilities.
\newblock {\em Statistics \& Probability Letters\/}~{\em 73\/}(2), 179 -- 187.

\bibitem[\protect\citeauthoryear{James}{James}{1991}]{James_1991}
James, S. (1991).
\newblock Approximate multinormal probabilities applied to correlated multiple
  endpoints in clinical trials.
\newblock {\em Stat Med\/}~{\em 10\/}(7), 1123--35.

\bibitem[\protect\citeauthoryear{Jasiulewicz and Kordecki}{Jasiulewicz and
  Kordecki}{2003}]{JK_2003}
Jasiulewicz, H. and W.~Kordecki (2003).
\newblock {Convolutions of Erlang and of Pascal distributions with applications
  to reliability }.
\newblock {\em Demonstratio Mathematica\/}~{\em 36}, 231--238.

\bibitem[\protect\citeauthoryear{Kost and McDermott}{Kost and
  McDermott}{2002}]{KM_2002}
Kost, J.~T. and M.~P. McDermott (2002).
\newblock Combining dependent p-values.
\newblock {\em Statistics \& Probability Letters\/}~{\em 60\/}(2), 183--190.

\bibitem[\protect\citeauthoryear{Koziol and Perlman}{Koziol and
  Perlman}{1978}]{KP_1978}
Koziol, J. and M.~Perlman (1978).
\newblock {Combining independent chi-squared tests.}
\newblock {\em J. Amer. Statist. Assoc.\/}~{\em 73}, 753--763.

\bibitem[\protect\citeauthoryear{Koziol}{Koziol}{1996}]{Koziol_1996}
Koziol, J.~A. (1996).
\newblock {A note on Lancaster's Procedure for the Combination of Independent
  Events}.
\newblock {\em Biometrical Journal\/}~{\em 38}, 653--660.

\bibitem[\protect\citeauthoryear{Koziol and Tuckwell}{Koziol and
  Tuckwell}{1994}]{KT_1994}
Koziol, J.~A. and H.~C. Tuckwell (1994).
\newblock {A weighted nonparametric procedure for the combination of
  independent events}.
\newblock {\em Biom. J.\/}~{\em 36}, 1005--1012.

\bibitem[\protect\citeauthoryear{Lancaster}{Lancaster}{1961}]{Lancaster_1961}
Lancaster, H.~D. (1961).
\newblock The combination of probabilities: an application of orthogonal
  functions.
\newblock {\em Austr. J. Statist.\/}~{\em 3}, 20--33.

\bibitem[\protect\citeauthoryear{Likes}{Likes}{1967}]{Likes_1967}
Likes, J. (1967).
\newblock Distributions of some statistics in samples from exponential and
  power-function populations.
\newblock {\em Journal of the American Statistical Association\/}~{\em
  62\/}(317), pp. 259--271.

\bibitem[\protect\citeauthoryear{Liptak}{Liptak}{1958}]{Liptak_1958}
Liptak, P. (1958).
\newblock {On the combination of independent tests}.
\newblock {\em Magyar Tud. Akad. Nat. Kutato int. Kozl.\/}~{\em 3}, 171--197.

\bibitem[\protect\citeauthoryear{Loesgen, Dempfle, Golla, and
  Bickeboller}{Loesgen et~al.}{2001}]{LDGB_2001}
Loesgen, S., A.~Dempfle, A.~Golla, and H.~Bickeboller (2001).
\newblock Weighting schemes in pooled linkage analysis.
\newblock {\em Genet Epidemiol\/}~{\em 21 Suppl 1\/}(NIL), S142--7.

\bibitem[\protect\citeauthoryear{MacNeill}{MacNeill}{1974}]{MacNeill_1974}
MacNeill, I.~B. (1974).
\newblock Tests for change of parameter at unknown times and distributions of
  some related functionals on brownian motion.
\newblock {\em The Annals of Statistics\/}~{\em 2\/}(5), pp. 950--962.

\bibitem[\protect\citeauthoryear{Mathai}{Mathai}{1983}]{Mathai_1983}
Mathai, A. (1983).
\newblock On linear combinations of independent exponential variables.
\newblock {\em Communications in Statistics - Theory and Methods\/}~{\em 12},
  625--632.

\bibitem[\protect\citeauthoryear{McGill and Gibbon}{McGill and
  Gibbon}{1965}]{MG_1965}
McGill, W.~J. and J.~Gibbon (1965).
\newblock The general-gamma distribution and reaction times.
\newblock {\em Journal of Mathematical Psychology\/}~{\em 2\/}(1), 1 -- 18.

\bibitem[\protect\citeauthoryear{Morse}{Morse}{1958}]{Morse_1958}
Morse, P. (1958).
\newblock {\em Queues, Inventories and Maintenance}.
\newblock John Wiley \& Sons.

\bibitem[\protect\citeauthoryear{Mosteller and Bush}{Mosteller and
  Bush}{1954}]{MB_1954}
Mosteller, F. and R.~R. Bush (1954).
\newblock {\em Selected quantitative techniques. In: Handbook of Social
  Psychology.}, Volume~1.
\newblock Cambridge, Mass.: Addison-Wesley.

\bibitem[\protect\citeauthoryear{Olkin}{Olkin}{1995}]{Olkin_1995}
Olkin, I. (1995).
\newblock Statistical and theoretical considerations in meta-analysis.
\newblock {\em Journal of Clinical Epidemiology\/}~{\em 48\/}(1), 133 -- 146.
\newblock The Potsdam International Consultation on Meta-Analysis.

\bibitem[\protect\citeauthoryear{Olkin and Saner}{Olkin and
  Saner}{2001}]{OS_2001}
Olkin, I. and H.~Saner (2001).
\newblock Approximations for trimmed {F}isher procedures in research synthesis.
\newblock {\em Stat Methods Med Res\/}~{\em 10\/}(4), 267--76.

\bibitem[\protect\citeauthoryear{Pearson}{Pearson}{1938}]{Pearson_1938}
Pearson, E.~S. (1938).
\newblock The probability integral transformation for testing goodness of fit
  and combining independent tests of significance.
\newblock {\em Biometrika\/}~{\em 30}, 134--148.

\bibitem[\protect\citeauthoryear{Pearson}{Pearson}{1933}]{Pearson_1933}
Pearson, K. (1933).
\newblock On a method of determining whether a sample of size n supposed to
  have been drawn from a parent population having a known probability integral
  has probably been drawn at random.
\newblock {\em Biometrika\/}~{\em 25}, 379--410.

\bibitem[\protect\citeauthoryear{Pepe and Fleming}{Pepe and
  Fleming}{1989}]{PF_1989}
Pepe, M.~S. and T.~R. Fleming (1989).
\newblock Weighted {K}aplan-{M}eier statistics: a class of distance tests for
  censored survival data.
\newblock {\em Biometrics\/}~{\em 45\/}(2), 497--507.

\bibitem[\protect\citeauthoryear{Pocock, Geller, and Tsiatis}{Pocock
  et~al.}{1987}]{PGT_1987}
Pocock, S.~J., N.~L. Geller, and A.~A. Tsiatis (1987).
\newblock The analysis of multiple endpoints in clinical trials.
\newblock {\em Biometrics\/}~{\em 43\/}(3), pp. 487--498.

\bibitem[\protect\citeauthoryear{Solomon and Stephens}{Solomon and
  Stephens}{1977}]{HS_1977}
Solomon, H. and M.~A. Stephens (1977).
\newblock Distribution of a sum of weighted chi-square variables.
\newblock {\em Journal of the American Statistical Association\/}~{\em
  72\/}(360), pp. 881--885.

\bibitem[\protect\citeauthoryear{Stouffer, Suchman, DeVinney, Star, and
  Williams}{Stouffer et~al.}{1949}]{Stouffer_1949}
Stouffer, S., E.~Suchman, L.~DeVinney, S.~Star, and R.~M.~J. Williams (1949).
\newblock {\em The American Soldier, Vol. 1: Adjustment during Army Life}.
\newblock Princeton: Princeton University Press.

\bibitem[\protect\citeauthoryear{Tippett}{Tippett}{1931}]{Tippett_1931}
Tippett, L. (1931).
\newblock {\em The Methods of Statistics}.
\newblock London: Williams and Norgate Ltd.

\bibitem[\protect\citeauthoryear{Wei and Lachin}{Wei and
  Lachin}{1984}]{WL_1984}
Wei, L.~J. and J.~M. Lachin (1984).
\newblock Two-sample asymptotically distribution-free tests for incomplete
  multivariate observations.
\newblock {\em Journal of the American Statistical Association\/}~{\em
  79\/}(387), pp. 653--661.

\bibitem[\protect\citeauthoryear{Whitlock}{Whitlock}{2005}]{Whitlock_2005}
Whitlock, M.~C. (2005).
\newblock Combining probability from independent tests: the weighted {Z}-method
  is superior to {F}isher's approach.
\newblock {\em J Evol Biol\/}~{\em 18\/}(5), 1368--73.

\bibitem[\protect\citeauthoryear{Yu, Gertz, Agarwala, Schaffer, and
  Altschul}{Yu et~al.}{2006}]{YGASetall_2006}
Yu, Y.-K., E.~M. Gertz, R.~Agarwala, A.~A. Schaffer, and S.~F. Altschul (2006).
\newblock Retrieval accuracy, statistical significance and compositional
  similarity in protein sequence database searches.
\newblock {\em Nucleic Acids Res\/}~{\em 34\/}(20), 5966--73.

\bibitem[\protect\citeauthoryear{Zelen and Joel}{Zelen and
  Joel}{1959}]{ZJ_1959}
Zelen, M. and L.~S. Joel (1959).
\newblock The weighted compounding of two independent significance tests.
\newblock {\em The Annals of Mathematical Statistics\/}~{\em 30\/}(4), pp.
  885--895.

\end{thebibliography}
\end{document}